\documentclass{amsart}

\usepackage{amsmath,amsfonts,amssymb,hyperref,mdwlist,graphicx,enumitem,color,multicol}
\usepackage[cmtip,all]{xy}
\usepackage[latin1]{inputenc}

\newtheorem{theorem}{Theorem}[section]
\newtheorem{lemma}[theorem]{Lemma}

\theoremstyle{definition}

\theoremstyle{remark}

\numberwithin{equation}{section}

\newcommand{\alg}{\textnormal{\bf{alg}}\,}
\newcommand{\vect}{\textnormal{\bf{span}}\,}
\newcommand{\am}{^{(+)}}\newcommand{\wt}{\widetilde}
\newcommand{\reg}{\mathcal{R}\textnormal{eg}\,}
\newcommand{\corf}{\mathbb{F}}\newcommand{\corfd}{\corf\cdot}
\newcommand{\jor}{\mathfrak{J}}\newcommand{\rad}{\mathcal{N}}\newcommand{\radd}{\mathcal{N}\,^2}
\newcommand{\josuu}{\jor\textnormal{osp}(1|2)}
\newcommand{\bim}{\mathcal{M}}\newcommand{\obim}{\mathcal{M}^{\textnormal{op}}}

\newcommand{\skewosp}{\mathcal{S}\textnormal{kew}(\alMndm,\osp)}
\newcommand{\skewospud}{\mathcal{S}\textnormal{kew}(\mathcal{M}_{1|2}(\corf),\osp)}

\newcommand{\heraa}{\mathcal{H}(\mathcal{A},\ast)}
\newcommand{\gra}{\Gamma}
\newcommand{\alga}{\mathcal{A}}
\newcommand{\algs}{\mathcal{S}}
 \newcommand{\mj}{\mathfrak{M}(\jor)}

\newcommand{\alMndm}{\mathcal{M}_{n|2m}(\corf)}\newcommand{\alMnm}{\mathcal{M}_{n|m}(\corf)}
\newcommand{\alMnn}{\mathcal{Q}(n)}
\newcommand{\osp}{\textnormal{osp}}
\newcommand{\algdt}{\mathcal{D}_{t}}

\newcommand{\josnm}{\jor\osp_{n|2m}(\corf)}
\newcommand{\josud}{\jor\osp_{1|2}(\corf)}

\newcommand{\suMnm}{\mathcal{M}_{n|m}(\corf)\am}

\newcommand{\nn}{e^{n}}\newcommand{\mm}{e^{m}}\newcommand{\dmdm}{e^{2m}}
\newcommand{\nm}{e^{nm}}\newcommand{\mn}{e^{mn}}
\newcommand{\ndm}{e^{n2m}}\newcommand{\dmn}{e^{2mn}}
\newcommand{\mdm}{e^{m2m}}\newcommand{\dmm}{e^{2mm}}
\newcommand{\hn}{h}\newcommand{\vmdm}{v}\newcommand{\smdm}{s}\newcommand{\sdmm}{\wt{s}}
\newcommand{\unm}{u}\newcommand{\kndm}{k}
\newcommand{\Hn}{H}\newcommand{\Vmdm}{V}\newcommand{\Smdm}{S}\newcommand{\Sdmm}{\wt{S}}
\newcommand{\Unm}{U}\newcommand{\Kndm}{K}

\newcommand{\bUnm}{\bar{U}}\newcommand{\bKndm}{\bar{K}}
\newcommand{\wUnm}{\wt{U}}\newcommand{\wKndm}{\wt{K}}

\newcommand{\gn}{\,g}\newcommand{\wmdm}{\,w}\newcommand{\zmdm}{\,z}\newcommand{\zdmm}{\,\wt{z}}
\newcommand{\ynm}{\,y}\newcommand{\xndm}{\,x}

\newcommand{\uu}{\eta^{u}}\newcommand{\kk}{\eta^{k}}

\newcommand{\ukg}{\eta^{uk,g}}\newcommand{\ukw}{\eta^{uk,w}}

\newcommand{\uk}{\eta^{uk}}\newcommand{\us}{\eta^{us}}
\newcommand{\ks}{\eta^{k\wt{s}}}
\newcommand{\uhy}{\eta^{uh,y}}\newcommand{\uhx}{\eta^{uh,x}}
\newcommand{\uvy}{\eta^{uv,y}}\newcommand{\uvx}{\eta^{uv,x}}
\newcommand{\usy}{\eta^{us,y}}\newcommand{\usx}{\eta^{us,x}}

\newcommand{\khy}{\eta^{kh,y}}\newcommand{\khx}{\eta^{kh,x}}
\newcommand{\kvy}{\eta^{kv,y}}\newcommand{\kvx}{\eta^{kv,x}}
\newcommand{\ksy}{\eta^{k\wt{s},y}}\newcommand{\ksx}{\eta^{k\wt{s},x}}

\begin{document}

\title{Orthosymplectic Jordan superalgebras and the Wedderburn principal theorem (WPT)}
\author{G\'omez  Gonz\'alez F.A. \& Vel\'asquez  R.}
\address[Faber G\'omez G.]{Calle 67 No. 53-108 Bloque 6-7 of.
 337}
 \email{faber.gomez@udea.edu.co}{}
\address[Raul Vel\'asquez]{Calle 67 No. 53-108 Bloque 6-7 of. 337}
\email{raul.velasquez@udea.edu.co}{}

\thanks{The authors were partially supported by Universidad de Antioquia, 
and the firts author by the proyect CODI 2014-1032}



\date{November 2016}

\maketitle

\begin{abstract}
An analogue of the Wedderburn Principal Theorem (WPT) 
is considered for finite-dimensional Jordan superalgebras 
$\alga$ with solvable radical 
$\rad,\, \radd=0$, 
and such that 
$\alga/\rad\cong\josnm$, 
where $\corf$ is a field of characteristic zero. 

\noindent
Let's we prove that the WPT is valid under some restrictions over the irreducible 
$\josnm$-bimodules contained in 
$\rad$, 
and it is shown with counter-examples that these restrictions can not be weakened.  
\end{abstract}

\subjclass{2001}{17C70,17C27,17C55}


\section*{Introduction}
In recent works, see  \cite{yo1} and \cite{yo2}, the first author proved an analogue to the 
Wedderburn principal theorem for Jordan superalgebras when we have a 
finite dimensional Jordan superalgebra 
$\alga$ with solvable radical $\rad$ such that
 $\radd=0$ and $\alga/\rad$ is a simple Jordan superalgebra of some of the following types: 
 superalgebra of superform,   Kac $\mathcal{K}_{10}$,  $\algdt$, $\suMnm$.
 Some conditions were impossed over the solvable radical $\rad$. 

Similarly as \cite{yo2}, we consider a finite dimensional Jordan superalgebra 
 $\alga$ 
over an algebraically closed field of characteristic 0 $\corf$, 
with solvable radical $\rad$ such that 
$\radd=0$ and $\alga/\rad\cong\josnm$, to follow we show that if 
$\rad$ considered as $\josnm$-superbimodule
does not contains any homomorphic image isomorphic to 
subbimodule $\reg(\josud)$ then, the Wedderburn principal theorem hold.  
Moreover, we shown that there is a counter-example to WPT for this case.

\section{Preliminary results and notations }
 Recall that an algebra $\alga$ is said to be a superalgebra if 
 $\alga=\alga_0\dotplus\alga_1$ satisfies the relation 
 $\alga_{i}\alga_{j}\subseteq \alga_{i+j(\text{\rm mod }2)}$, i.e.
  $\alga$ is a $\mathbb{Z}_2$ - graded algebra.  
  Given an element $a\in\alga_0\cup\alga_1$, $|a|=i$ denotes its parity, 
  according to $a\in\alga_i$.

\medskip Let 
$\Gamma=\alg\langle\,1,e_i,\,\, i\in\mathbb{Z}^{+} | e_ie_j+e_je_i=0\,\rangle$ 
denote the Grassmann algebra. 
Then, $ \gra=\gra_0\dotplus\gra_1$, where $\gra_0$ and $\gra_1$, 
is spanned by all monomials of even and odd length respectively, 
and it is easy to see that  
$\gra$ has a superalgebra structure.

\medskip Let $\alga=\alga_0\dotplus\alga_1$ be a superalgebra, let's the Grassmann envelope of 
$\alga$ is the algebra 
$\gra(\alga)=\gra_0\otimes\alga_0+\gra_1\otimes\alga_1$.  
Assume that 
$\mathfrak{M}$ is a homogeneous variety of algebras (see, eg. \cite{4rusos}).  
The superalgebra $\alga$  is said to be a  
$\mathfrak{M}$-superalgebra if the Grassmann envelope 
$\gra(\alga)$ belongs to $\mathfrak{M}$. 

\medskip 
An associative superalgebra is just a $\mathbb{Z}_2$-graded associative algebra, 
but it is not so in general (see \cite{Shes}). 
One can easily check that a superalgebra  
$\jor=\jor_0\dotplus\jor_1$
 is a Jordan superalgebra if and only if it satisfies the super identities 
\begin{align}
&\label{idsupercom} xy=(-1)^{|x||y|}yx,\\
&\label{idsupjordan}((xy)z)t+(-1)^{|t|(|z|+|y|)+|z||y|}((xt)z)y+(-1)^{|x|(|y|+|z|+|t|)+|t||z|}((yt)z)x\\
&\quad=(xy)(zt)+(-1)^{|t||z|+|t||y|}(xt)(yz)+(-1)^{|y||z|}(xz)(yt).\nonumber
\end{align} 
In particular, the Jordan superalgebra 
$\jor=\jor_0\dotplus\jor_1$ is a ($\mathbb{Z}_2$-graded) 
Jordan algebra if and only if 
$(\jor_1)^2=0$.
\medskip\noindent 

\medskip 
Let 
$\alga$ 
be an associative superalgebra with multiplication 
$ab$, we consider a new multiplication 
$a\circ b=\frac{1}{2}(ab+(-1)^{|a||b|}ba)$, $a,\,b\in\alga_0\cup\alga_1$. 
We can see, that with respect to this multiplication 
$\alga$ 
has a structure of Jordan superalgebra, which we will denote as 
$\alga\am$. 
 
\par \medskip \noindent C.T.C. Wall proved in \cite{Wal} that every associative simple finite-dimensional
 superalgebra over an algebraically closed field 
 $\corf$ is isomorphic to one of following associative superalgebras 
\begin{itemize}
\item[{\rm (i)}] 
$\alga=\alMnm$,\quad $\alga_0=\Biggl\{\begin{bmatrix} a& 0\\ 0 & d \end{bmatrix}\Biggr\},\quad
 \alga_1=\Biggl\{\begin{bmatrix} 0& b\\ c & 0\end{bmatrix}\Biggr\},$ 
\item[{\rm (ii)}] 
$\alga=\alMnn $,\quad $\alga_0=\Biggl\{\begin{bmatrix} a& 0\\ 0 & a \end{bmatrix}\Biggr\},\quad
 \alga_1=\Biggl\{\begin{bmatrix}0& h\\ h & 0\end{bmatrix}\Biggr\}.$
\end{itemize}
where 
$a, \, h\in\mathcal{M}_{n}(\corf)$, $d\in\mathcal{M}_{m}(\corf)$, 
$b\in\mathcal{M}_{n\times m}(\corf)$, $c\in\mathcal{M}_{m\times n}(\corf)$. 

\medskip Let 
$\alga$ be an associative (super)algebra.
 A superinvolution $\ast$ in $\alga$ 
 is a graded linear mapping 
 $\ast:\alga\longrightarrow\alga$ such that 
 $(a^\ast)^\ast=a$ and $(ab)^{\ast}=(-1)^{|a||b|}b^{\ast}a^{\ast}$. 
 Let $\heraa$ be the set of symmetric elements of $\alga$ relative to $\ast$, 
 namely, $\heraa=\{a\in\alga/ a^\ast=a\}$.
 It is easy to see that $\heraa\subseteq\alga\am$, and therefore 
 $\heraa$ is a Jordan superalgebra.
 
Let $I_n$, $I_m$ be identity matrices of order $n$ and $m$ respectively.  
We denote by $t$ the usual transposition of matrices and let 
$$U=-U^t=-U^{-1}=\begin{bmatrix} 0 & -I_m\\ I_m & 0\end{bmatrix} .$$

Consider the linear mapping $\osp:\alMndm\longrightarrow \alMndm$, given by
\begin{equation*}
\begin{bmatrix}a& b \\ c& d \end{bmatrix}^\osp =
\begin{bmatrix} I_n & 0\\ 0& U \end{bmatrix}\begin{bmatrix} a^t & -c^t\\ b^t & d^t \end{bmatrix} \begin{bmatrix}
I_n& 0\\ 0& U^{-1}\end{bmatrix},
\end{equation*} where $a\in\bim_{n}$, $b, c^t\in\bim_{n\times 2m}$ and $d\in\bim_{2m}$.

\medskip We can see that 
$\osp$ is a superinvolution over superalgebra $\alMndm$. 
So, the Jordan superalgebra 
$\mathcal{H}(\alMndm,\osp)$, denote by 
$\josnm$, is determined by the following  matrices set:

\begin{equation*}
\josnm=\begin{Bmatrix} & 
\begin{bmatrix}a& b_1& b_2 \\ -b_2^t & d_1 & d_2 \\ 
b_1^t& d_3& d_1^t \end{bmatrix} & \Big|&    a=a^t,\, d_2=-d_2^t, \, d_3=-d_3^t & 
\end{Bmatrix}
\end{equation*} 
where 
$a\in\mathcal{M}_n(\corf)$, 
$b_1$, $b_2\in\bim_{n\times m}(\corf)$, 
$d_1$, $d_2,$ 
and 
$d_3\in\bim_m(\corf)$.

\medskip 
Simple finite-dimensional Jordan superalgebras over zero characteristic fields 
were classified by V.G. Kac \cite{Kac} (see also \cite{Kan}).

Now, we recall that a 
$\jor$-superbimodule $\bim=\bim_0\dotplus\bim_1$ 
is a Jordan superbimodule if the corresponding split null extension 
$\mathcal{E}=\jor\oplus\bim$ is a Jordan superalgebra. 
Besides, the split null extension is the vector space direct sum 
$\jor\oplus\bim$ 
with multiplication that extends the multiplication of 
$\jor$, the action of $\jor$ on $\bim$, and $\bim^2=0$.  
Let $\bim$ be a 
$\jor$-superbimodule, 
the opposite superbimodule $\obim=\obim_0\dotplus\obim_1$ 
is defined by the conditions $\obim_0=\bim_1,\,\obim_1=\bim_0$, 
and the following action of 
$\jor$, 
$a\cdot m^{\textnormal{op}}=(-1)^{|a|}(am)^{\textnormal{op}}, \, 
m^{\textnormal{op}}\cdot a=(ma)^{\textnormal{op}}$,
for any 
$a\in\jor_0\cup\jor_1,\,m\in\obim_0\cup\obim_1$. 
Whenever $\bim$ is a Jordan 
$\jor$-superbimodule it is possible to see that 
$\obim$ is so as well. A {\it{regular superbimodule}} 
$\reg(\jor)$ is defined on the vector superspace 
$\jor$ with the action of 
$\jor$ coinciding with the multiplication in 
$\jor$.

\medskip 
Irreducible bimodules over Jordan superalgebra $\josnm$ 
were classified by C. Martinez and E. Zelmanov in \cite{zelmar2}, 
who proved that the only unital irreducible 
$\josnm$-bimodules are the regular bimodule 
$\reg(\josnm)$  the bimodule 
$\mathfrak{S}=\mathcal{S}\textnormal{kew}(\alMndm,\osp)$, and their opposites, 
where
\begin{equation}\label{skewgral}
\mathfrak{S}=
\begin{Bmatrix} & 
\begin{bmatrix} 
a& b_1& b_2 \\ b_2^t & d_1 & d_2 \\ -b_1^t& d_3& -d_1^t 
\end{bmatrix} & \Big|&    a=-a^t,\, d_2=d_2^t, \, d_3=d_3^t & 
\end{Bmatrix}
\end{equation} 

$a\in\mathcal{M}_n(\corf)$, $b_1$, $b_2\in\bim_{n\times m}(\corf)$,  $d_1$, $d_2,$ 
and $d_3\in\bim_m(\corf)$.

\medskip
\subsection{The Peirce decomposition} 
Recall, that if $\jor$ is a Jordan (super)algebra with unity $1$, and  $\{e_1,\ldots,e_n\}$ is a 
set of pairwise orthogonal idempotents such that $1=\sum_{i=1}^n e_i$, then $\jor$ 
admits Peirce decomposition \cite{zelmar2},  
it is
$$\jor=\biggl( \bigoplus_{i=1}^n\jor_{ii} \biggr) \bigoplus \biggl( \bigoplus_{i<j}\jor_{ij}\biggr) , $$ 
where 
$$\jor_{ii}=\{\,x\in\jor: \quad e_ix=x,\}$$ and 
$$\jor_{ij}=\{\,x\in\jor: \quad e_ix=\frac{1}{2}x,\quad e_jx=\frac{1}{2}x\,\}, \quad i\neq j$$ 
are the Peirce components of $\jor$ relative to the idempotents $e_i$, and $ e_j$, 
moreover the following relations hold
\begin{align*}
&\jor_{ij}^2\subseteq \jor_{ii}+\jor_{jj},\quad \jor_{ij}\cdot\jor_{jk}\subseteq \jor_{ik}\\ 
&\jor_{ij}\cdot\jor_{kl}=0\quad \textnormal{when}\quad i\neq k,l; j\neq k,l.
\end{align*}

\section{Main theorem}
In this section we prove the central result of this paper. 
To start we introduce the following notation, by 
$e_{ij}$, $i,j=1,\ldots , n+2m$, we denote the usual unit matrices. 

\noindent For $i,j\in\{1,\ldots, n\}$ and $p,q\in\{1,\ldots,m\}$, denote  
$\nn_{ij}=e_{ij}$, $\nm_{ip}=e_{i\,n+p}$, $\ndm_{i\,p}=e_{i\,n+m+p}$,
$\mm_{pq}=e_{n+p\,n+q}$, $\mn_{pi}=e_{n+p\,i}$,  $\mdm_{pq}=e_{n+p\,n+m+q}$,
$\dmdm_{pq}=e_{n+m+p\,n+m+q}$, $\dmm_{pq}=e_{n+m+p\,n+q}$ 
and $\dmn_{pi}=e_{n+m+p\,i}$.

\medskip \noindent 
Consider 
$\hn_{ij}=\nn_{ij}+\nn_{ji}$ if $i\neq j$, 
$\hn_{ii}=\nn_{ii}$, 
$\vmdm_{pq}=\mm_{pq}+\dmdm_{qp}$,
$\smdm_{pq}=\mdm_{pq}-\mdm_{qp}$, 
$\sdmm_{pq}=\dmm_{pq}-\dmm_{qp}$, 
$\unm_{ip}=\nm_{ip}+\dmn_{pi}$,
$\kndm_{ip}=\ndm_{ip}-\mn_{pi}$.

\medskip With the previous notation,
the Jordan superalgebra
$\jor=\josnm$ is spanned by 
$\{\hn_{ij}, \vmdm_{pq}, \smdm_{pq}, \sdmm_{pq}, \unm_{ip},\,\kndm_{ip}\}$ 
and its dimension is given by
$\frac{(n+2m)^2+n-2m}{2}$. 

From 
$a\circ b=\frac{1}{2}(ab+(-1)^{|a||b|}ba)$ 
we can see that the non-zero products of basis elements of  
$\jor$ are defined  as follows:

\scalebox{0.86}{\parbox[c][3cm][l]{1.1\linewidth}{ %
\begin{align}\label{mjosp0}
\hn_{ii}\circ \hn_{ii}=\hn_{ii}, \quad \hn_{ij} \circ \hn_{kl}=\frac{1}{2}(\delta_{jk}\hn_{il}+\delta_{li}\hn_{kj}+
&\delta_{jl}\hn_{ik}+\delta_{ik}\hn_{jl}), \textnormal{ if } i\neq j, k\neq l,\nonumber \\
 \smdm_{pq}\circ \sdmm_{rt}=\dfrac{1}{2}(\delta_{qr}\vmdm_{pt}+\delta_{pt}\vmdm_{qr}-
 \delta_{qt}\vmdm_{pr}-\delta_{pr}\vmdm_{qt}), \quad &
  \vmdm_{pq}\circ\vmdm_{rt}=\frac{1}{2}(\delta_{qr}\vmdm_{pt}+\delta_{pt}\vmdm_{rq}),\nonumber\\
\vmdm_{pq}\circ \smdm_{rt}=\frac{1}{2}(\delta_{qr}\smdm_{pt}+\delta_{tq}\smdm_{rp}), \quad
 \vmdm_{pq}\circ \sdmm_{rt}=&\frac{1}{2}(\delta_{pr}\sdmm_{qt}+\delta_{pt}\sdmm_{rq}),
\end{align}}}

\scalebox{0.86}{\parbox[c][3cm][l]{1.1\linewidth}{ %
\begin{align}\label{mjosp01}
\unm_{kr}\circ \hn_{ij}=\dfrac{1}{2}(\delta_{jk}\unm_{ir}+\delta_{ik}\unm_{jr}), \quad 
 \kndm_{lp}\circ\hn_{ij}=\frac{1}{2}(\delta_{jl}\kndm_{ip}+ &
\delta_{il}\kndm_{jp}),  \textnormal{ if } i\neq j \nonumber\\
\unm_{kr}\circ \hn_{ii}=\dfrac{1}{2}\delta_{ik}\unm_{ir},\,  \kndm_{lp}\circ\hn_{ii}=\frac{1}{2}\delta_{il}\kndm_{ip},\,
\unm_{kr}\circ \vmdm_{pq}=\frac{1}{2}\delta_{rp}&\unm_{kq}, \,  \kndm_{lr}\circ
 \vmdm_{pq}=\frac{1}{2}\delta_{rq}\kndm_{lp},\nonumber\\
 \unm_{ir}\circ\smdm_{pq}=\frac{1}{2}(\delta_{rp}\kndm_{iq}-\delta_{rq}\kndm_{ip}), \quad 
  \kndm_{ir}\circ\sdmm_{pq}=\frac{1}{2}(\delta_{rp}\unm_{iq}&-\delta_{rq}\unm_{ip}), 
\end{align}}}

\scalebox{0.86}{\parbox[c][3cm][l]{1.1\linewidth}{ %
\begin{align}\label{mjosp1}
\unm_{ip}\circ \unm_{jq}=\frac{1}{2}\delta_{ij}\sdmm_{pq},\quad  \kndm_{ip}\circ \kndm_{jq}=\frac{1}{2}\delta_{ij}&
 \smdm_{qp}, \quad \unm_{ip}\circ\kndm_{iq}=\frac{1}{2}\vmdm_{qp}-\delta_{pq}\hn_{ii},\nonumber \\
\unm_{ip}\circ\kndm_{jq}=&\frac{1}{2}(\delta_{ij}\vmdm_{qp}-\delta_{pq}\hn_{ij}) \textnormal{ if } i\neq j, 
\end{align}}}

where $\delta_{ij}=0$ if $i\neq j$ and $\delta_{ii}=1$. 
We note that  
the products in \eqref{mjosp0} and \eqref{mjosp01} are symmetric
 and the products in \eqref{mjosp1} is skew-symmetric.

Now we prove the following theorem

\emph{\bfseries{Theorem}}
Let $\alga$ be a finite-dimensional Jordan superalgebra with a solvable radical 
$\rad$, such that  $\radd=0$, and 
the quotient superalgebra 
$\jor=\alga/\rad$ is isomorphic ti $\josnm$. 
If $n+m\geq 3$ or $n=m=1$ and  no homomorphic image of $\rad$, 
considered as a $\jor$-bimodule, contains a subbimodule isomorphic to 
$\reg(\josud)$, then there exists a subsuperalgebra $\algs\subset\alga$ such that 
$\algs\cong\josnm$ and $\alga=\algs\oplus\rad$.

\proof 
Take $\jor=\josnm$ and let $\jor$-mod denote the category of Jordan $\jor$-superbimodules. 
By Theorem 8.1 in \cite{zelmar2}, every $V\in \jor$-mod is completely reducible. 
Let $\mj$ by the set of $V$ in $ \jor$-mod such that  
$V$ does not contain a bimodule isomorphic to
$\reg(\josud)$,\ among its irreducible summands. 
Clearly, $\mj$ is closed with respect to subbimodules and homomorphic images, 
and by \cite{yo1} (Theorem 3.3)  
we observe that it is suffices to prove the theorem when 
$\alga$ is unital and $\rad$ is irreducible. 
Following  \cite{zelmar2}, (Theorem 6.3), there are four different
types of unital irreducible $\jor$-bimodules 
$reg(\josnm)$, $\text{Skew}(\mathcal{M}_{n\mid 2m}, \osp)$ 
and their opposites.

We have that 
$$\josnm=
\textnormal{alg}\langle\, \hn_{ij},\, \vmdm_{pq},\, \smdm_{pq},\, \sdmm_{pq}\,\rangle
\dotplus 
\textnormal{vect}\langle\, \unm_{ip},\, \kndm_{ip}\,\rangle,$$
where $i,j=1,\ldots, n$ and $p,q=1,\ldots, m$.

Since WPT is valid for Jordan algebras, 
and using the fact that 
$\alga/\rad\cong\josnm$ 
we can assume that there exists 
$\algs_0\subseteq\alga_0$ such that 
$\algs_0\cong\alg\langle\,\hn_{ij}, \,\vmdm_{pq},\, \smdm_{pq}, \, \sdmm_{pq}\,\rangle$.
Therefore there exist 
$\Hn_{ij}, \,\Vmdm_{pq},\, \Smdm_{pq}$, and $\Sdmm_{pq}\in\alga_0$ 
for which the multiplications \eqref{mjosp0} are valid when we substitute
$\hn_{ij}, \,\vmdm_{pq}$, $\smdm_{pq}$ and 
$\sdmm_{pq}$ by $\Hn_{ij}, \,\Vmdm_{pq},\, \Smdm_{pq}$ 
and $\Sdmm_{pq}$,  
respectively.

\medskip We note that 
$\{\Hn_{ii}, \Vmdm_{pp}$  for $ i=1\ldots, n; p=1\ldots, m\}$
 is a set of pairwise orthogonal idempotents such that 
 $\Hn_{11}+\cdots+\Hn_{nn}+\Vmdm_{11}+\cdots+\Vmdm_{mm}=1$.
 Thus $\alga$ has a Peirce decomposition with respect  to itsidempotents:
 
$$\alga=\Biggl( \bigoplus_{\substack{i\leq j \\ i=1}}^n
(\alga)_{ij}\Biggr)\bigoplus\Biggl(\bigoplus_{\substack{i=1\ldots n \\
 p=1\ldots, m}} (\alga)_{pq}\Biggr)
 \bigoplus\Biggl( \bigoplus_{\substack{p\leq q\\ p=1}}^m(\alga)_{pq}\Biggr)\,.$$

Now we need to find 
$\wUnm_{ip}$ and $\wKndm_{ip}\in\alga_1$ 
such that the multiplications \eqref{mjosp01} and \eqref{mjosp1} 
hold when we change 
$\hn_{ij}, \,\vmdm_{pq}$, $\smdm_{pq}$, $\sdmm_{pq}$, $\unm_{ip}$ and  
$\kndm_{ip}$ by $\Hn_{ij}$, $\Vmdm_{pq}$, $\Smdm_{pq}$, $\Sdmm_{pq}$, $\wUnm_{ip}$ 
and $\wKndm_{ip}$, respectively.

Since $\alga_1/\rad_1\cong(\josnm)_1$,
there exist $\bUnm_{ip}$ and $\bKndm_{ip}\in\alga_1/\rad_1$  
such that \eqref{mjosp01} and \eqref{mjosp1} are valid in 
$\alga/\rad$ when we change  
$\hn_{ij}, \,\vmdm_{pq}$, $\smdm_{pq}$, $\sdmm_{pq}$, $\unm_{ip}$ and  
$\kndm_{ip}$ by $\Hn_{ij}, \,\Vmdm_{pq},\, \Smdm_{pq}$, $\Sdmm_{pq}$, $\bUnm_{ip}$ 
and $\bKndm_{ip}$,  respectively.

\medskip\textbf{Case 1}  $\rad\cong\reg(\josnm)$. 

\medskip
Let  $\gn_{ij}$, $\wmdm_{pq}$, $\zmdm_{pq}$, $\zdmm_{pq}$, $\ynm_{ip}$ 
and 
$\xndm_{ip}\in\rad$ and assume that the isomorphism $\sigma$ is determinated by 

\scalebox{0.85}{\parbox[c][1.3cm][l]{1.1\linewidth}{ %
\begin{equation*}
\begin{aligned}
 \sigma(\gn_{ij})= \hn_{ij}, \, \sigma(\wmdm_{pq} )=\vmdm_{pq}, \, 
 \sigma (\zmdm_{pq}) = \smdm_{pq}, \\ 
 \sigma(\zdmm_{pq}) =\sdmm_{pq}, \,
\sigma( \ynm_{ip})=\unm_{ip},  \sigma(\xndm_{ip})= \kndm_{ip}.
\end{aligned}
\end{equation*}}} 

So, we have that

\scalebox{0.86} {\parbox[c][1.3cm][l]{1.1\linewidth}{ %
$$\rad_0=\vect\langle
\gn_{ij},\, \wmdm_{pq},\, \zmdm_{pq},\, \zdmm_{pq}\rangle \textnormal{ and } \rad_1=\vect\langle 
 \ynm_{ip},\,\xndm_{ip}\rangle.$$}}

  The action of 
  $\josnm$ over $\rad$ is determined by the equations 
  \eqref{mjosp0},\eqref{mjosp01} and \eqref{mjosp1} when we replace 
  $\gn_{ij}$, $\wmdm_{pq}$, $\zmdm_{pq}$, $\zdmm_{pq}$, $\ynm_{ip}$,
   and $\xndm_{ip}$ by 
   $\hn_{ij}$, $\vmdm_{pq}$, $\smdm_{pq}$, $\sdmm_{pq}$, $\unm_{ip}$, and $\kndm_{ip}$, 
   respectively.

\begin{lemma}\label{lemmac11}
 Let $\varphi:\alga\longrightarrow \alga/\rad$ be  the canonical homomorphism.
For $i=1, \cdots, n$ and $p=1,\cdots, m$
let 
$U_{ip}$ and $K_{ip}\in\alga_1$ 
be preimages of 
$\bar{U}_{ip}$ and $\bar{K}_{ip}$ 
respectively, then 

\scalebox{0.86}{\parbox[c][1.3cm][l]{1.1\linewidth}{ %
\begin{align}\label{uhkv}
\Unm_{ip}\cdot \Hn_{ij}=
\frac{1}{2}\Unm_{jp},\quad  \Unm_{ip}\cdot \Vmdm_{pq}=\frac{1}{2}\Unm_{iq}, \, &
\Kndm_{ip}\cdot \Hn_{ij}=
\frac{1}{2}\Kndm_{jp},\quad  \Kndm_{ip}\cdot \Vmdm_{qp}=\frac{1}{2}\Kndm_{iq} 
\end{align} }}

\scalebox{0.86}{\parbox[c][1.3cm][l]{1.1\linewidth}{ %
\begin{align} \label{rela1}
\Unm_{ip}\cdot \Smdm_{pq}=
\frac{1}{2}\Kndm_{iq}, \quad \Kndm_{ip}\cdot \Sdmm_{pq}=\frac{1}{2}\Unm_{iq}
\end{align}}}
\end{lemma} 

\proof 
To start we prove \eqref{uhkv}. 
From 
 $\varphi(\Unm_{ip}\cdot \Hn_{ij})=\frac{1}{2}\bUnm_{jp}$, 
 $\varphi(\Unm_{ip}\cdot \Vmdm_{pq})=\varphi(\Unm_{iq}\cdot \Vmdm_{qq})=\frac{1}{2}\bUnm_{iq}$ 
 and using the properties of Peirce decomposition for the Jordan superalgebra 
 $\alga$, we note that 
 $\Unm_{ip}\cdot \Hn_{ij}\in(\alga_0)_{j\,n+p}$. 
 We can see that $\{\ynm_{jp},\xndm_{jp}\}$ is a generator set of $(\rad_0)_{jp}$, 
 and therefore we can assume that there exist 
 $\uhx_{ipij},\,\uhy_{ipij}\in\corf$ such that 
\begin{align}
\Unm_{ip}\cdot \Hn_{ij}=\frac{1}{2}\Unm_{jp}+\uhy_{ipij}\ynm_{jp}+\uhx_{ipij}\xndm_{jp}. \label{uiphij}\end{align}

Similarly, 
 \begin{align}
K_{ip}\cdot \Vmdm_{qp}=\frac{1}{2}K_{iq}+\kvy_{ipqp}\ynm_{iq}+\kvx_{ipqp}\xndm_{iq},\label{uipvpq}
\end{align} for some $\uvx_{ippq},\,\uvy_{ippq}\in\corf$.

Using \eqref{mjosp0} and replacing 
$x=\Unm_{ip}$, $y=z=t=\Hn_{ii}$ in \eqref{idsupjordan} we have 
\begin{equation}\label{L11}
2((\Unm_{ip}\cdot \Hn_{ii})\cdot \Hn_{ii})\cdot \Hn_{ii}+\Hn_{ii}\cdot \Unm_{ip}=3(\Unm_{ip}\cdot
 \Hn_{ii})\cdot\Hn_{ii}.\end{equation}
  If we replace \eqref{uiphij} in \eqref{L11}, 
 we obtain  $$\frac{3}{4}\Unm_{ip}+\frac{5}{2}\uhy_{ipii}\ynm_{ip}+
 \frac{5}{2}\uhx_{ipii}\xndm_{ip}=\frac{3}{4}\Unm_{ip}+3\uhy_{ipii}\ynm_{ip}+3\uhx_{ipii}\xndm_{ip}.$$ 
 Hence, 
 $\uhx_{ipii}\xndm_{ip}+\uhy_{ipii}\ynm_{ip}=0$. 
 Since $\ynm_{ip}$  and $\xndm_{ip}$ are linearly independent, we have  
 $\uhy_{ipii}=\uhx_{ipii}=0$. Similarly, we can prove that $\khy_{ipii}=\khx_{ipii}=0$. 

\medskip \noindent Using the fact that   
$\varphi(\Unm_{ip}\cdot \Hn_{ij})=\varphi(\Unm_{jp}\cdot \Hn_{jj})$  
we can conclude that   
$\uhx_{ipij}=\uhx_{jpjj}=$ $\uhy_{ipij}=\uhy_{jpjj}=0$.

\medskip Now, if we replace 
$x=\Kndm_{ip}$, $y=z=t=\Vmdm_{pp}$ in \eqref{idsupjordan} and using \eqref{mjosp0} we obtain
 $$2((\Kndm_{ip}\cdot \Vmdm_{pp})\cdot \Vmdm_{pp})\cdot
  \Vmdm_{pp}+\Vmdm_{pp}\cdot\Kndm_{ip}=3(\Kndm_{ip}\cdot \Vmdm_{pp})\cdot\Vmdm_{pp}.$$ 
As in the case above, it is easy to see that 
 $\kvy_{ippp}=\kvx_{ippp}=0$.
  Again, using some properties of the canonical homomorphism we obtain 
 $\kvy_{ipqp}=\kvx_{ipqp}=\kvx_{ippp}=0$.  
 The other equalities in \eqref{uhkv} are proven similarly.

\noindent\medspace Now, we show the equality \eqref{rela1}. 

Let $\usy_{ippq}$, $\usx_{ippq}$, $\ksx_{ippq}$ and $\ksy_{ippq}$ scalars such that 

\scalebox{0.86}{\parbox[c][1.3cm][l]{1.1\linewidth}{ %
\begin{equation}\label{usks}
\Unm_{ip}\cdot \Smdm_{pq}=\frac{1}{2}\Kndm_{iq}+\usy_{ippq}y_{iq}+\usx_{ippq}x_{iq}, \quad \Kndm_{ip}\cdot
 \Sdmm_{pq}=\frac{1}{2}\Unm_{iq}+\ksy_{ippq}y_{iq}+\ksx_{ippq}x_{iq}.
\end{equation}}}

Let $x=\Unm_{ip}$, $y=\Smdm_{pq}$, $z=\Vmdm_{pq}$ and $t=\Sdmm_{pq}$ 
in the equation \eqref{idsupjordan}. 
Using  \eqref{mjosp0} and \eqref{uhkv}  we have, 
\begin{align}\label{usukigual}
((\Unm_{ip}\cdot \Smdm_{pq})\cdot \Vmdm_{pq})\cdot \Sdmm_{pq}-\frac{1}{4}\Unm_{iq}=-\frac{1}{8}\Unm_{iq}
\end{align}
Using \eqref{usks}  and  computing the products we obtain 
 $(\ksy_{iqpq}-\usx_{ippq})\ynm_{iq}=0$  and $\ksx_{ippq}x_{iq}=0$. 
  Hence $\ksx_{ippq}=0$ (Similarly, we can prove 
  $\usy_{ippq}=0$), and 
 \begin{align}\label{igualdadusks}
 \ksy_{iqpq}-\usx_{ippq}=0. 
\end{align} 

\noindent\medspace We note that
 $\varphi(\Unm_{ir}\cdot\Smdm_{rq})=
 \varphi(\Unm_{ip}\cdot\Smdm_{pq})=-\varphi(\Unm_{ip}\cdot\Smdm_{qp})=\frac{1}{2}\bKndm_{iq}$. 
 Thus we have $\usx_{irrq}=\usx_{ippq}=-\usx_{ipqp}$; 
 but this relation only depend of $i$ and $q$ so we can write 
\begin{align}\label{usq}
\Unm_{ip}\cdot \Smdm_{pq}=\frac{1}{2}\Kndm_{iq}+\us_{iq}x_{iq}.
\end{align}  
Analogously we have 
$\Kndm_{ip}\cdot \Sdmm_{pq}=\frac{1}{2}\Unm_{iq}-\us_{iq}y_{iq}$. 
The equation \eqref{igualdadusks} can be rewritten as $ \ks_{ip}=\us_{iq}$, 
only depend of $i$ and therefore $ \ks_{i}=\us_{i}$. 
Since 
$K_{ip}\cdot \Sdmm_{pq}=-K_{ip}\cdot \Sdmm_{qp}$,
 we obtain 
 $\ks_{i}=-\ks_{i}=0$ similarly, $\us_{i}=0$.  
 Hence, 
 $\Unm_{ip}\cdot \Smdm_{pq}=\frac{1}{2}\Kndm_{iq}$, $\Kndm_{ip}\cdot \Sdmm_{pq}=\frac{1}{2}\Unm_{iq}$.
 \endproof

\begin{lemma}\label{lemmac12}
 Let $\varphi:\alga\longrightarrow \alga/\rad$ be  the canonical homomorphism.
For $i=1, \cdots, n$ and $p=1,\cdots, m$
let 
$U_{ip}$ and $K_{ip}\in\alga_1$ 
be preimages of 
$\bar{U}_{ip}$ and $\bar{K}_{ip}$ 
respectively, then 

\scalebox{0.86}{\parbox[c][1.3cm][l]{1.1\linewidth}{ %
\begin{align} \label{rela2} 
\Unm_{ip}\cdot \Unm_{iq}=
\frac{1}{2}\Sdmm_{pq},   \quad  \Kndm_{jp}\cdot \Kndm_{jq}=\frac{1}{2}\Smdm_{qp}, 
\end{align}}}

\scalebox{0.86}{\parbox[c][1.3cm][l]{1.1\linewidth}{ %
\begin{align}\label{rela4}
\Unm_{ip}\cdot\Kndm_{iq}=
\frac{1}{2}\Vmdm_{qp},\quad \Unm_{ip}\cdot\Kndm_{jp}=-\frac{1}{2}\Hn_{ij}, \quad
\Unm_{ip}\cdot\Kndm_{ip}=
\frac{1}{2}\Vmdm_{pp}-\Hn_{ii}.
\end{align} }} 
\end{lemma}

\medskip Now we prove \eqref{rela2}. Let $\eta^{u,\wt{z}}_{ipq}$ and  $\eta^{u,z}_{ipq}\in\corf$,  
such that $\Unm_{ip}\cdot
 \Unm_{iq}=\frac{1}{2}\Sdmm_{pq}+\eta^{u,\wt{z}}_{ipq}\zdmm_{pq}+\eta^{u,z}_{ipq}\zmdm_{pq}.$
We note that 
$\varphi(\Unm_{ip}\cdot \Unm_{iq})=\varphi(\Unm_{jp}\cdot \Unm_{jq})$ 
then 
 $\eta^{u,\wt{z}}_{ipq}$ only depend  of 
 $p$ and $q$
 therefore $\eta^{u,\wt{z}}_{jpq}=\eta^{u,\wt{z}}_{pq}$ for all  $j=1,\ldots, n$.  
  Now  we can write 
$\Unm_{ip}\cdot \Unm_{iq}=\frac{1}{2}\Sdmm_{pq}+\eta^{u,\wt{z}}_{pq}\zdmm_{pq}+\eta^{u,z}_{pq}\zmdm_{pq}.$

Substituting $x$, $y$, $z$ and $t$ in the equation \eqref{idsupjordan} respectively by 
$\Unm_{ip},$ $\Unm_{iq}$, $\Sdmm_{qp}$ and $\Vmdm_{pq}$, 
and using \eqref{mjosp0} we obtain   
 $\Unm_{iq}\cdot\Vmdm_{pq}=\Unm_{iq}\cdot\Sdmm_{qp}=\Vmdm_{pq}\cdot\Sdmm_{rs}=0$,
 so because of \eqref{uhkv}   we obtain  
$((\Unm_{ip}\cdot\Unm_{iq})\cdot \Sdmm_{qp})\cdot \Vmdm_{pq}=0.$  
Hence $\eta^{u,z}_{pq} \wmdm_{pq}=0$ and therefore $\eta^{u,{z}}_{pq}=0$. So, we have 
\begin{equation}\label{uu}
\Unm_{ip}\cdot \Unm_{iq}=\frac{1}{2}\Sdmm_{pq}+\eta_{pq}^{u,\wt{z}}\zdmm_{pq}.
\end{equation}

By a similar process we can prove that 
\begin{equation}\label{kk}
\Kndm_{ip}\cdot \Kndm_{iq}=\frac{1}{2}\Smdm_{qp}+\eta_{pq}^{k,z}\zmdm_{qp}.
\end{equation}

\medskip Now we can consider the product $\Unm_{ip}\cdot\Kndm_{jp}$. 
Let $\ukg_{ijp}\in\corf$ such that  
$\Unm_{ip}\cdot\Kndm_{jp}=-\frac{1}{2}\Hn_{ij}+\ukg_{ijp}g_{ij}.$  
Knowing that 
 $\varphi(\Unm_{ip}\cdot\Kndm_{jp})=\varphi(\Unm_{jp}\cdot\Kndm_{ip})=\varphi(\Unm_{iq}\cdot\Kndm_{jq})$ 
we  can affirm that   $\ukg_{ijp}=\ukg_{jip}=\ukg_{ijq}.$ 
Similarly, we have $\uk_{ipq}=\uk_{jpq}$. So we can write  
\begin{align}\label{ukuk}
\Unm_{ip}\cdot\Kndm_{jp}=-\frac{1}{2}\Hn_{ij}+\uk_g g_{ij} \textnormal{ and }
 \Unm_{ip}\cdot\Kndm_{iq}=\frac{1}{2}\Vmdm_{qp}+\uk_{pq} w_{qp}.
\end{align}

If we replace $x,\,y,\, z, $ and $t$, respectively, by 
$\Unm_{ip},$ $\Unm_{iq}$, $\Smdm_{qp}$ and 
$\Vmdm_{pq}$ in the equation \eqref{idsupjordan}, and using \eqref{mjosp0} and \eqref{uhkv} we obtain 
$$((\Unm_{ip}\cdot\Unm_{iq})\cdot \Smdm_{qp})\cdot \Vmdm_{pq}+\frac{1}{2}(\Unm_{iq}\cdot \Smdm_{qp})\cdot
 \Unm_{iq}=\frac{1}{2}\Unm_{iq}\cdot (\Unm_{iq}\cdot \Smdm_{qp})$$ 
Replacing \eqref{usq} and \eqref{uu} in the equality above we obtain 
 $\uu_{pq}-\uk_{qp}=0$. Similarly we can show that  
 $0=\kk_{pq}-\uk_{pq}.$ 
\begin{align}\label{uusuk}
\uu_{pq}=\uk_{qp}, \textnormal{ and } \kk_{pq}=\uk_{pq}
\end{align}

\noindent Using \eqref{uusuk} we obtain $\uu_{pq}+\kk_{pq}=\uk_{qp}+\uk_{pq}.$ 
We note that  $\uu_{pq}=-\uu_{qp}$ and $\kk_{pq}=-\kk_{qp}$. 
Thus $\uu_{qp}+\kk_{qp}=\uk_{pq}+\uk_{qp}=\uu_{pq}+\kk_{pq},$ 
hence  $\uu_{pq}+\kk_{pq}=0$ and therefore
 $\uk_{pq}+\uk_{qp}=0$.

\medspace Finally, we show the equality \eqref{rela4}. 
Let $x=\Unm_{ip}$, $y=\Kndm_{ip}$ and $z=t=\Hn_{ij}$ 
in the equality \eqref{idsupjordan}. As we did before, we have  
$$((\Unm_{ip}\cdot\Kndm_{ip})\cdot \Hn_{ij})\cdot \Hn_{ij} + \frac{1}{2}\Unm_{ip}\cdot\Kndm_{ip} =
(\Unm_{ip}\cdot \Kndm_{ip})\cdot {\Hn_{ij}}^2+\frac{1}{2}\Unm_{jp}\cdot\Kndm_{jp} \label{ukip}.$$ 
Thus 
$(\ukw_{jp}-\ukw_{ip})w_{pp}+(\ukg_{jp}-\ukg_{ip})g_{jj}=0$. 
Since $w_{pp}$ and $g_{ii}$ 
are linearly independent we have 
$\ukw_{jp}=\ukw_{ip}=\ukw_p$ and $\ukg_{jp}=\ukg_{ip}=\ukg_p$.

If we take 
$x=\Unm_{ip}$, $y=\Kndm_{ip}$,  $z=\Vmdm_{pq}$ and 
$t=\Vmdm_{qp}$ in \eqref{idsupjordan} we can show that  
$\ukw_{p}=\ukw_{q}=\uk_w$ and $\ukg_{p}=\ukg_{q}=\uk_g$ and therefore, we obtain 
\begin{align}\label{baseuukk}
\Unm_{ip}\cdot\Kndm_{ip}=\frac{1}{2}\Vmdm_{pp}-\Hn_{ii}+\uk_w w_{pp}+\uk_g g_{ii}.
\end{align}
 Let $x=\Unm_{ip}$, $y=\Unm_{iq}$, $z=\Kndm_{ip}$ and $t=\Kndm_{iq}$ 
 in the equation \eqref{idsupjordan}. Thus we obtain 

\scalebox{0.86}{\parbox[c][2.5cm][l]{1.1\linewidth}{ %
\begin{align}\label{uukk}
((\Unm_{ip}\cdot&\Unm_{iq})\cdot\Kndm_{ip})\cdot\Kndm_{iq}+\Unm_{iq}\cdot((\Unm_{ip}\cdot\Kndm_{iq})\cdot
 \Kndm_{ip})+((\Unm_{iq}\cdot\Kndm_{iq})\cdot\Kndm_{ip})\cdot\Unm_{ip}=\\
 (\Unm_{ip}\cdot&\Unm_{iq})\cdot(\Kndm_{ip}\cdot\Kndm_{iq})+(\Unm_{iq}\cdot\Kndm_{ip})\cdot
  (\Unm_{ip}\cdot\Kndm_{iq})+(\Unm_{iq}\cdot\Kndm_{iq})\cdot(\Kndm_{ip}\cdot\Unm_{ip})\nonumber
\end{align}}} 

If we replace \eqref{uu}, \eqref{kk} and \eqref{baseuukk} in \eqref{uukk} and using the fact 
$\us_{p}=\ks_{p}=0$ and $\uu_{pq}=\uk_{qp}=-\uk_{pq}$ and $\kk_{pq}=\uk_{pq}$ we obtain 
$$(2\uk_w-2\uu_{pq})w_{qq}+(2\uk_w-\uk_g)w_{pp}-2\uk_g g_{ii}=0.$$ 
Due to the fact that $w_{pp}$, $w_{qq}$ and $g_{ii}$ are linearly independent, 
$\uu_{pq}=\uk_w=\uk_g=0$. Hence $\uk_{pq}=\kk_{pq}=0$. \endproof

We note that the relations \eqref{uhkv} - \eqref{rela4} are valid when 
$n=1$ and $m\geq 2$ or $n\geq 2$ and $m=1$ and therefore, the WPT is valid.

\medskip
\textbf{Case 2 } $\rad\cong\skewosp$. 

We denote 
\begin{equation*}
\begin{aligned}
a_{ij}=\nn_{ij}-\nn_{ji},\quad \wt{a}_{pq}=\mm_{pq}-\dmdm_{qp},\quad f_{pq}=\mdm_{pq}+\mdm_{qp}, \\
 \wt{f}_{pq}=\dmm_{pq}+\dmm_{qp},\quad b_{ip}=\nm_{ip}-\dmn_{pi},\quad c_{ip}=\ndm_{ip}+\mn_{pi}.
\end{aligned}
\end{equation*}

Since \eqref{skewgral} we can see that $\rad$ is generated by 
$a_{ij}$, $\wt{a}_{pq}$, $f_{pq}$, $\wt{f}_{pq}$, $b_{ip}$ and $c_{ip}$, 
moreover, 
$$\rad_0=\vect\langle a_{ij}, \wt{a}_{pq}, f_{pq}, \wt{f}_{pq}\rangle, \text{ and } 
\rad_1=\vect\langle b_{ip}, \, c_{ip}\rangle$$
 
It is easy to see that the action of 
$\josnm$ over $\skewosp$ 
is determined by the following multiplication table

\scalebox{0.86}{\parbox[c][5cm][l]{1.1\linewidth}{ %
\begin{align}\label{skew0}
\hn_{ii}\circ a_{ii}=\frac{1}{2}a_{ii}, \quad \hn_{kl} \circ
 a_{ij}=\frac{1}{2}(\delta_{jk}a_{il}+\delta_{li}a_{kj}+\delta_{jl}a_{ik}+\delta_{ik}a_{lj}), 
\textnormal{ if } i\neq j, k\neq l,\nonumber \\
 \sdmm_{pq}\circ
  f_{rt}=\frac{1}{2}(\delta_{pr}\wt{a}_{tq}+\delta_{ps}\wt{a}_{rq}-\delta_{qr}\wt{a}_{tp}-\delta_{qt}\wt{a}_{pr}), \quad
  \sdmm_{pq}\circ \wt{a}_{rt}=\frac{1}{2}( \delta_{qr}\wt{f}_{pt}-\delta_{pr}\wt{f}_{qt}),  
 \nonumber \\
\smdm_{pq}\circ
 \wt{f}_{rt}=\dfrac{1}{2}(\delta_{qr}\wt{a}_{pt}+\delta_{qt}\wt{a}_{pr}-\delta_{pr}\wt{a}_{qt}-\delta_{pt}\wt{a}_{qr}), \quad 
\smdm_{pq}\circ \wt{a}_{rt}=\frac{1}{2}(\delta_{pt}f_{qr}-\delta_{qt}f_{pr}), \nonumber \\
\vmdm_{pq}\circ\wt{a}_{rt}=\frac{1}{2}(\delta_{qr}\wt{a}_{pt}+\delta_{pt}\wt{a}_{rq}), \quad  \vmdm_{pq}\circ
 f_{rt}=\frac{1}{2}(\delta_{qr}f_{pt}+\delta_{tq}f_{pr}), \nonumber\\  
\vmdm_{pq}\circ  \wt{f}_{rt}=\frac{1}{2}(\delta_{pr}\wt{f}_{qt}+\delta_{pt}\wt{f}_{qr}), \quad 
\end{align} }}

\scalebox{0.86}{\parbox[c][6.5cm][l]{1.1\linewidth}{ %
\begin{align}\label{skew01}
h_{ij}\circ b_{kr}=\dfrac{1}{2}(\delta_{jk}b_{ir}+\delta_{ik}b_{jr}), \quad h_{ii}\circ b_{kr}=\dfrac{1}{2}\delta_{ik}b_{ir},
 \quad  \vmdm_{pq}\circ b_{kr}=\frac{1}{2}\delta_{rp}b_{kq},\nonumber \\ 
h_{ij}\circ c_{kr}=\frac{1}{2}(\delta_{jk}c_{ir}+\delta_{ik}\kndm_{jr}), \quad  h_{ii}\circ
 c_{kr}=\frac{1}{2}\delta_{ik}c_{ir},\quad   v_{pq}\circ c_{kr}=\frac{1}{2}\delta_{rq}c_{kp},\nonumber \\
s_{pq}\circ b_{ir}=\frac{1}{2}(\delta_{pr}c_{iq}-\delta_{qr}c_{ip}),\quad \wt{s}_{pq}\circ
 c_{ir}=\frac{1}{2}(\delta_{rp}b_{iq}-\delta_{qr}b_{ip})\nonumber \\
\unm_{ip}\circ a_{kj}=\frac{1}{2}(\delta_{ij}b_{kp}-\delta_{ik}b_{jp}),\quad \unm_{ip}\circ
 \wt{a}_{qr}=\frac{1}{2}\delta_{pq}b_{ir}, \nonumber \\
\kndm_{ip}\circ a_{jk}=\frac{1}{2}(\delta_{ik}c_{jp}-\delta_{ij}c_{kp}),\quad \kndm_{ip}\circ
 \wt{a}_{qr}=-\frac{1}{2}\delta_{pr}c_{iq},\quad \nonumber \\
\unm_{ip}\circ f_{qr}=\frac{1}{2}(\delta_{pq}c_{ir}+\delta_{pr}c_{iq}),\quad \kndm_{ip}\circ
 \wt{f}_{qr}=\frac{1}{2}(\delta_{pq}b_{ir}+\delta_{pr}b_{iq}), \nonumber\\ 
 \unm_{ir}\circ\smdm_{pq}=\frac{1}{2}(\delta_{rp}\kndm_{iq}-\delta_{rq}\kndm_{ip}), \quad 
  \kndm_{ir}\circ\sdmm_{pq}=\frac{1}{2}(\delta_{rp}\unm_{iq}-\delta_{rq}\unm_{ip}), 
\end{align}}}

\scalebox{0.86}{\parbox[c][3cm][l]{1.1\linewidth}{ %
\begin{align}\label{skew1}
\unm_{ip}\circ b_{jq}=\frac{1}{2}\delta_{ij}\wt{f}_{pq},\quad  u_{ip}\circ 
 c_{jq}=\frac{1}{2}(\delta_{pq}a_{ij}-\delta_{ij}\wt{a}_{qp})\nonumber \\
\kndm_{ip}\circ c_{jq}=-\frac{1}{2}\delta_{ij}f_{pq} , \quad  k_{ip}\circ 
 b_{jq}=\frac{1}{2}(\delta_{pq}a_{ji}-\delta_{ij}\wt{a}_{pq}).
\end{align}}}

The products in \eqref{skew0}, \eqref{skew01} are symmetric and the products in 
\eqref{skew1} are skew-symmetric. 
 
Similarly as was prove to case 1, we have that if $n+m\geq 3$ 
then there exist an analogous to Lemma \ref{lemmac11} and \ref{lemmac12}. Therefore
there exist $U_{ip}$ and $K_{ip}\in\alga_1$ such that 
$\rm{vect}\langle U_{ip}, K_{ip}, \, i=1,\ldots, n,\, p=1,\ldots, m\rangle\cong(\josnm)_1$, 
thus the WPT is valid.

\noindent Now we prove the theorem for $m=n=1$ and $\rad\cong\skewospud$. 
 Let $h$, $v$, $\bar{u}$ and $\bar{k}$ such that 
 $(\corfd h+\corfd v)\dotplus (\corfd \bar{u}+\corfd \bar{k)} \cong\josud$. 
 We need to find $\wt{u}$, $\wt{k}\in\alga_1$ 
 such that $\varphi({\wt{u}})=\bar{u}$ and $\varphi(\wt{k})=\bar{k}$. 
 Moroever, $\wt{u}h=\wt{u}v=\frac{1}{2}u$, 
$\wt{k}h=\wt{k}v=\frac{1}{2}k$,  and $\wt{u}\wt{k}=\frac{1}{2}v-h$.

Let $\rad=(\corfd \wt{a}+\corfd f+\corfd \wt{f})\dotplus(\corfd b+\corfd c)$ 
and consider the following action of $\josud$ over $\rad$
\begin{align}\label{josuuskewpar}
v\wt{a}=\wt{a}, \, vf=f, \, v\wt{f}=\wt{f}, \quad & bh=bv=\frac{1}{2}b, \quad  ch=cv=\frac{1}{2}c, \nonumber \\ 
\wt{u}\wt{a}=\wt{k}\wt{f}=\frac{1}{2}b, \quad & \wt{k}\wt{a}=\wt{u}f=\frac{1}{2}c   
\end{align}
\begin{equation}\label{josuuskewimpar}
\wt{u}b=\wt{f}, \quad \wt{u}c=-\frac{1}{2}\wt{a} = \wt{k}b, \quad  \wt{k}c=-f, 
\end{equation} 
where \eqref{josuuskewpar} and \eqref{josuuskewimpar} 
are commutative and anticommutative, respectively. 
It is easy to see that $\rad\cong\skewospud$.

\medskip
Let $\xi_{\wt{a}}$, $\xi_{f}$ and $\xi_{\wt{f}}\in\corf$ such that 
 $uk=\frac{1}{2}v-h+\xi_{\wt{a}}\wt{a}+\xi_{f}f+\xi_{\wt{f}}\wt{f}$. 
 Let's prove that there exist 
 $\alpha_b$, $\alpha_c$, $\beta_b$ and $\beta_c\in\corf$ such that 
$$\wt{u}=u+\alpha_b b+\alpha_c c,\quad \wt{k}=k+\beta_b b+\beta_c c,\text{ and }
\wt{u}\cdot\wt{k}=\frac{1}{2}v-h.$$
 
 We note that 
 $\varphi(\wt{u})=\bar{u}$ and $\varphi(\wt{k})=\bar{k}$. 
 Using  \eqref{josuuskewpar} we have 
 $$\wt{u}h=\wt{u}v=\frac{1}{2}\wt{u} \text{ and } \wt{k}h=\wt{k}v=\frac{1}{2}\wt{k}.$$

 Now $\wt{u}\wt{k}=\frac{1}{2}v-h$ if and only if 
$$ (\xi_{\wt{a}}+\frac{1}{2}\alpha_b-\frac{1}{2}\beta_c)\wt{a}+(\xi_{f}+\alpha_c)f+(\xi_{\wt{f}}+\beta_b)\wt{f}=0.$$
Since $\wt{a}$, $f$, and $\wt{f}$ are linearly independent we have
 $$\xi_{\wt{a}}+\frac{1}{2}\alpha_b-\frac{1}{2}\beta_c=\xi_{f}+\alpha_c=\xi_{\wt{f}}+\beta_b=0,$$ 
 and therefore 
 $2\xi_{\wt{a}}+\alpha_b=\beta_c$, $\xi_{f}=-\alpha_c$, $\xi_{\wt{f}}=-\beta_b$ 
 is a solution, hence the WPT is valid.

We note that if $\rad$ is isomorphic to anyone of superbimodules opposites,  
then by the Pierce properties, we have that the radical part in any product is zero 
and therefore  the equalities
\eqref{uhkv}, \eqref{rela1}, \eqref{rela2} and \eqref{rela4} hold 
when we change  
$H_{ij}$, $V_{pq}$, $S_{pq}$, $\wt{S}_{pq}$, $U_{ij}$ and $K_{ip}$
respectively by
${h}_{ij}$, ${v}_{pq}$, $s_{pq}$, $\wt{s}_{pq}$, $u_{ij}$ and $k_{ip}$, 
and therefore, WPT is true.

\section{counter-example} 

Let $\alga=\alga_0\oplus\alga_1$ where $\alga_0=\corfd h+\corfd v+\corfd g+\corfd w$ and 
$\alga_1=\corfd u+\corfd k+\corfd y+\corfd x$, and  
 $\rad=(\corfd g+\corfd w)\dotplus(\corfd y+\corfd x)$. 
 The non-zero multiplications in $\alga$ are given by 

\scalebox{0.86}{\parbox[c][1.5cm][l]{1.1\linewidth}{ %
\begin{equation}\label{counter1_110}
\begin{array}{llll}
h^2=h, \, v^2=v, \, &  hg=g, \, vw=w, \, & uh=uv=\frac{1}{2}u,  &  kh=k v= \frac{1}{2}k,\\
yh=yv=\frac{1}{2}y, \, & ug=uw=\frac{1}{2}y,\, & xh=xv= \frac{1}{2}x, \, & kg=kw=\frac{1}{2}x, 
\end{array}
\end{equation}}}

\scalebox{0.86}{\parbox[c][0,9cm][l]{1.1\linewidth}{ %
\begin{equation}\label{counter1_111}
\begin{array}{lll}
ux=\frac{1}{2}w-g, \quad & yk=\frac{1}{2}w-g, \quad & u k=\frac{1}{2}v-h+g ,
\end{array}
\end{equation}}}

\noindent 
where the products \eqref{counter1_110} are symmetric 
and \eqref{counter1_111} are skew-symmetric.

\medskip 
\noindent 
Using \eqref{idsupjordan}, and the table of multiplications above
it is easy to show that $\alga$ is a Jordan superalgebra. 
Moreover, if
$\josud_0=\corfd h_{11}+\corfd v_{11}$ while $\josud_{1}=\corfd u_{11}+\corfd k_{11}$.
Consider the mapping
$\varphi:\alga/\rad \longrightarrow \josud$ 
and 
$\psi:\rad\longrightarrow \josud$ 
given by 
$\varphi(h)=\psi(g)=h_{11}$, 
$\varphi(v)=\psi(w)=v_{11}$, 
$\varphi(u)=\psi(y)=u_{11}$ 
and $\varphi(k)=\psi(x)=k_{11}$.

\noindent 
We can see that $\varphi$ is an isomorphism between $\alga/\rad$ and $\josud$, 
while $psi$ is an isomorphism between $\rad$ and $\reg(\josuu)$.

If we assume that the WPT is valid for $\alga$,  
then there exist ${h}$, ${v}\in\alga_0$ and $\wt{u}$, $\wt{k} \in\alga_1$ 
such that, the following products are commutative
${h}^2={h}$, ${v}^2=v$, ${h}\wt{u}={v}\wt{u}=\frac{1}{2}\wt{u}$, 
${h}\wt{k}={v}\wt{k}=\frac{1}{2}\wt{k}$, 
and anticommutative product
$\wt{u}\wt{k}=\frac{1}{2}\wt{v}-\wt{h}$ 
hold, 
and $\wt{u}\equiv u (\textrm{mod } \rad)$ and 
$\wt{k}\equiv k (\textrm{mod } \rad)$. 

\medskip 
Consider $\alpha_x,\,\alpha_y, \, \beta_x$ and $\beta_y\in\corf$ 
such that 
$\wt{u}=u+\alpha_y y+\alpha_x x$ and 
$\wt{k}=k+\beta_y y+\beta_x x$.  We note that 
\begin{align*}
\wt{u}\wt{k}=(u+\alpha_y y+\alpha_xx)& (k+\beta_x x+\beta_yy)
 = uk +\alpha_y y k+\beta_x u x\\ 
 =&\frac{1}{2}v-h+g+\frac{1}{2}\alpha_y w-\alpha_y g+\frac{1}{2}\beta_x w-\beta_x g\\
 =& \frac{1}{2}v-h+(1-\alpha_y-\beta_x)g+\frac{1}{2}(\alpha_y+\beta_x)w. 
\end{align*} 
So $\wt{u}\wt{k}=\frac{1}{2}v-h$ if and only if $2(1-\alpha_y-\beta_x)g+(\alpha_y+\beta_x)w=0$. 
Due to $g$ and $w$ are linearly independent, we have $1=\alpha_y+\beta_x=0$ and so we have a contradiction.

\end{document}